\documentclass[10pt,reqno]{amsart}
\usepackage{amssymb}
\usepackage{amsfonts}
\usepackage{amsthm}
\usepackage{amsmath}
\usepackage{hyperref}

\newcommand{\Reals}{{\mathbb{R}}}
\newcommand{\Cmplx}{{\mathbb{C}}}
\newcommand{\Ints}{{\mathbb{Z}}}
\newcommand{\Disk}{{\mathbb{D}}}

\let\Im=\undefined\DeclareMathOperator*{\Im}{Im}

\DeclareMathOperator*{\supp}{supp} 
\DeclareMathOperator{\Exp}{{{\mathbb{E}}}}

\newtheorem{theorem}{Theorem}[section]
\newtheorem{prop}[theorem]{Proposition}
\newtheorem{lemma}[theorem]{Lemma}

\theoremstyle{definition}
\newtheorem{definition}[theorem]{Definition}
\theoremstyle{remark}

\numberwithin{equation}{section}

\begin{document}

\title[Gaussian fluctuations for $\beta$~Ensembles.]{Gaussian fluctuations for $\beta$~Ensembles.}
\author[R.~Killip]{Rowan Killip}
\address{Rowan Killip\\
         UCLA Mathematics Department\\
         Box 951555\\
         Los Angeles, CA 90095}
\email{killip@math.ucla.edu}

\date{\today}

\begin{abstract}
We study the Circular and Jacobi $\beta$-Ensembles and prove Gaussian fluctuations for the number of points
in one or more intervals in the macroscopic scaling limit.
\end{abstract}

\maketitle

\section{Introduction}

The circular $\beta$ ensemble with $n$ points, or C$\beta$E${}_n$, was introduced by Dyson, \cite{Dyson},
as a simplification of previous random matrix ensembles.  It is a random process of $n$ points distributed
on the unit circle in $\Cmplx$ so that 
\begin{align}\label{E:CbetaE}
\Exp_n^\beta (f) &= \frac{1}{Z_{n,\beta}} \int_{-\pi}^\pi\!\! \cdots  \int_{-\pi}^\pi
    f(e^{i\phi_1},\ldots,e^{i\phi_n})
    \bigl|\Delta(e^{i\phi_1},\ldots,e^{i\phi_n})\bigr|^\beta \, \frac{d\phi_1}{2\pi} \cdots \frac{d\phi_n}{2\pi}
\end{align}
for any symmetric function $f$. Here, $\Delta$ denotes the Vandermonde determinant,
\begin{equation}\label{VDefn}
  \Delta(z_1,\ldots,z_n) = \prod_{1\leq j < k \leq n} \!\! (z_k - z_j),
\end{equation}
and the partition function is given by
\begin{equation}\label{CGpart}
  Z_{n,\beta} = \frac{\Gamma(\tfrac12\beta n + 1)}{\bigl[\Gamma(\tfrac12\beta + 1)\bigr]^n};
\end{equation}
as shown in \cite{Good,Wilson}.  This law also arises as the Gibbs measure for $n$
identical charged particles confined to lie on the circle and interacting via the $2$-dimensional
Coulomb law.  For this reason, it is also known as the $\log$-gas.  An important example of \eqref{E:CbetaE}
is as the distribution of eigenvalues of an $n\times n$ unitary matrix chosen at random
according to Haar measure; this corresponds to $\beta=2$.

The second point process we will discuss is known as the Jacobi $\beta$-ensemble with $n$-points.  (We
will abbreviate this name to J$\beta$E$_n$.)  In this case, there are $n$-points on the interval $[-2,2]$
with joint probability density
\begin{align}\label{E:JbetaE}
p_n^\beta(x_1,\ldots,x_n) &= \frac{1}{Z} \bigl|\Delta(x_1,\ldots,x_n)\bigr|^\beta
        \prod_{j=1}^n \, (2-x_j)^{a-1}(2+x_j)^{b-1},
\end{align}
where the parameters $a,b$ are positive real numbers. 
The partition function, $Z$, was determined by Selberg, \cite{Selberg}.
The name Jacobi ensemble derives from the relation to the classical orthogonal polynomials of this name.

In this paper, we determine the statistical behaviour of the number of particles in one or more intervals
as $n\to\infty$ with the length of the intervals remaining fixed as $n\to\infty$.  This is what may
be termed \textit{macroscopic} statistics; the \textit{microscopic} regime, where the length of the intervals
scale as $\frac1n$, is much more difficult and will not be addressed.

The statement for the Jacobi ensemble is the cleaner of the two:

\begin{theorem}\label{T:mainJ}
Given a sample from J$\beta$E$_n$, let $N_n(\theta)$ denote the number of points that lie in the interval $[2\cos(\theta),2]$.
Then for any distinct $\theta_1,\ldots,\theta_J\in(0,\pi)$,
$$
\sqrt{\tfrac{\pi^2\beta}{\log(n)}}[N_n(\theta_j) - n\theta_j],\qquad j\in\{1,\ldots,J\}
$$
converge to independent Gaussian random variables with mean $0$ and variance $1$ as $n\to\infty$.
\end{theorem}

For the circular ensemble, we need one preliminary.  Let us write $\Psi$ for the zero-mean
Gaussian process on $(-\pi,\pi)$ with covariance
$$
\Exp\{\Psi(\theta_1)\Psi(\theta_2)\} = \begin{cases} 1 &: \theta_1=\theta_2 \\  0 &: \theta_1\neq\theta_2 \end{cases}
$$
that is, $\Psi(\theta)$ represents an independent Gaussian associated to each point in this arc.

\begin{theorem}\label{T:mainC}
Let us fix $\theta_1<\theta_2<\cdots<\theta_J\in(-\pi,\pi)$ and write $N(a,b)$ for the number of points
in a sample from C$\beta$E$_n$ that lie in the arc between $a$ and $b$.  Then the joint distribution of
$$
\sqrt{\tfrac{\pi^2\beta}{\log(n)}}\bigl[N_n(\theta_j,\theta_{l})- n\tfrac{\theta_l-\theta_j}{2\pi} \bigr], \qquad j<l
$$
converges to that of $\Psi(\theta_l) - \Psi(\theta_j)$ as $n\to\infty$.
\end{theorem}

In the special case of a single arc, we see that
\begin{equation}\label{OneArc}
\sqrt{\tfrac{\pi^2\beta}{2\log(n)}}\bigl[N_n(0,\theta)- \tfrac{n\theta}{2\pi} \bigr]
\end{equation}
converges to Gaussian random variable of mean zero and variance one.

For $\beta\in\{1,2,4\}$, Theorem~\ref{T:mainC} was proved by Costin and Lebowitz, \cite{CostLeb}.
At these three temperatures, the model becomes exactly soluble in the sense that there are
explicit determinantal formulae for the correlation functions; see \cite{ForresterBook} or \cite{Mehta}.
There are many other works on Gaussian fluctuations at these three temperatures; in particular, we
would like to draw the reader's attention to \cite{DiaEvans,DiaShah,Joh2,Jonsson,Sosh1,Sosh2,Sosh3,Wieand}.
One of the questions addressed in these papers are the laws of other linear statistics.  Given a point process
$\{x_j\}$, say on $\Reals$, and a function $f:\Reals\to\Reals$, the associated linear statistic is
$$
X_f := \sum_j f(x_j).
$$
In particular, the number of points in an interval $I$ corresponds to $f=\chi_I$.

For the circular problem with $\beta=2$, the behaviour of linear statistics can be recast as a question about
the asymptotics of Toeplitz determinants.  In this case, we may attribute the proof of asymptotically Gaussian fluctuations
to Szeg\H{o},~\cite{SzegoToe}, at least for $f\in C^{1+\epsilon}$.  Of the numerous papers devoted to Szeg\H{o}'s
Theorem in the last half-century, one stands out for its adaptability to the case of general $\beta$,
namely \cite{Joh1}.  This paper proves asymptotically Gaussian fluctuations for C$\beta$E with $f\in C^{1+\epsilon}$.
A subsequent paper, \cite{Joh2}, covers ensembles on the real line.

The approach taken here is to use matrix models for C$\beta$E${}_n$ and J$\beta$E${}_n$ that where described in
\cite{KN}.  This paper extended work of Dumitriu and Edelman, \cite{DumE},
who discovered tri-diagonal matrix models the $\beta$-Hermite and $\beta$-Laguerre ensembles.  By `matrix model' we mean
an ensemble of random matrices whose eigenvalues follow the desired law.  We should also mention \cite{Trotter}
as the progenitor of both these papers.

In \cite{DumE2}, Dumitriu and Edelman proved Gaussian linear statistics for polynomials by studying
traces of powers of their matrix models.

In \cite{KStoi}, the circular matrix models where studied in the microscopic scaling limit using a Pr\"ufer variables approach.
Having finished that paper, it dawned upon me that the use of Pr\"ufer variables leads to a very simple treatment of
the problem of the number of particles in an interval.  This is what is presented here.

The virtue of using Prufer variables to treat the microscopic scaling limit of the matrix models was discovered
independently by Valk\'o and Vir\'ag, \cite{VV}.  Specifically, they study the eigenvalue statistics of the
$\beta$-Hermite ensemble in a neighbourhood of a reference energy lying in the bulk.
Earlier, \cite{RRV} studied the statistics in a neighbourhood of the edge using the Riccati transformation. 

\subsection*{Notation}  While dealing with the Jacobi ensembles, we will write $X\lesssim Y$ to indicate that $X\leq C Y$
for some constant $C$. In all cases, the implicit constant will depend only on $\beta$, $a$, and $b$.

\subsection*{Acknowledgements}
The author was supported in part, by NSF grant DMS-0401277 and a Sloan Foundation Fellowship.   He is also
grateful to the Institute for Advanced Study (Princeton) for its hospitality.

\section{Background}

\subsection{A change of variables}
In this subsection, we recount the results of \cite{KN}.  This paper describes probability distributions
on discrete measures with the property that the marginal distribution of the location of the mass points
happen to have the laws that interest us.  The virtue of the random measures described in \cite{KN}
appears when one looks at them from the point of view of orthogonal polynomials.

Given a probability measure $d\mu$ on the unit circle in $\Cmplx$ that is supported at exactly $N$ points, we may write
\begin{equation}\label{E:CmeasDefn}
\int f\,d\mu = \sum  \mu_j f(e^{i\theta_j}).
\end{equation}
Applying the Gram--Schmidt procedure to $\{1,z,\ldots,z^{N-1}\}$ leads to an orthogonal basis for $L^2(d\mu)$ built of monic
polynomials.  We write $\Phi_0(z)\equiv 1,\ldots,\Phi_{N-1}(z)$ for these polynomials.

As discovered by Szeg\H{o}, these polynomials obey a recurrence relation,
\begin{equation}
\begin{aligned}\label{SzegoRec}
\Phi_{k+1}(z) &= z\Phi_{k}(z) - \bar\alpha_k \Phi_{k}^*(z) \\
\Phi_{k+1}^*(z) &= \Phi_{k}^*(z) - \alpha_k z\Phi_{k}(z)
\end{aligned}
\end{equation}
with $\Phi_0^*=\Phi_0=1$.  Here $\Phi_k^*$ denotes the reversed polynomial,
\begin{equation}\label{PhiReflect}
\Phi_{k}^*(z) = z^k \overline{\Phi_{k}(\bar z^{-1})},
\end{equation}
and $\alpha_0,\ldots,\alpha_{N-2}$ are recurrence coefficients belonging to $\Disk$.  Expanding
$z^N$ in this basis shows that
\begin{equation}\label{Phi_N}
z\Phi_{N-1}(z) - e^{-i\eta} \Phi^*_{N-1}(z) = 0 \qquad\text{in $L^2(d\mu)$}
\end{equation}
for a unique $\eta\in[0,2\pi)$.  A basic fact in the theory of orthogonal polynomials is that the mapping
$$
\{ (\theta_1,\mu_1),\ldots,(\theta_N,\mu_N) \} \mapsto (\alpha_0,\ldots,\alpha_{N-2},e^{-i\eta})  
$$
is a bijection. For a proof of this (and the other claims in this subsection),
see \cite{SimonOPUC1} or \cite{Szego}.  As in the former reference, we will refer to $\alpha_0,\ldots,\alpha_{N-2},e^{-i\eta}$
as the Verblunsky parameters of $d\mu$. 

It will be expedient for us to encode \eqref{SzegoRec} in a different way.  As the zeros of the orthogonal
polynomials lie inside the unit disk, $B_k(z) = z\Phi_k(z)/\Phi^*_k(z)$ is a Blaschke product of degree $k+1$; moreover,
\begin{equation}\label{BlaschkeRec}
B_0(z) = z,
\qquad
B_{k+1}(z) = z B_k(z) \frac{1-\bar\alpha_k \bar B_k(z)}{1-\alpha_k B_k(z)}.
\end{equation}
Lastly, by \eqref{Phi_N}, we have $\supp(d\mu) = \{ z : B_{N-1}(z)= e^{-i\eta}\}$.

We have two systems of coordinates for the set of probability measures on $S^1$ with exactly $N$ mass points, on the one hand,
the location and weight of the masses, on the other, the Verblunsky coefficients.  Choosing random Verblunsky coefficients
gives rise to a random measure and thus to a random subset of the unit circle.  The good news is that
they will be chosen independently.  In the model for the the circular ensembles, they will follow a $\Theta$ distribution:

\begin{definition}  A complex random variable, $X$, with values in the unit disk,
$\Disk$, is $\Theta_\nu$-distributed (for $\nu>1$) if
\begin{equation}\label{E:ThetaDefn}
\Exp\{f(X)\} = \tfrac{\nu-1}{2\pi} \int\!\!\!\int_\Disk f(z) (1-|z|^2)^{(\nu-3)/2} \,d^2z.
\end{equation}
Simple computations show $\Exp\{X\}=0$, $\Exp\{ |X|^2 \} = \tfrac{2}{\nu+1}$ and
$\Exp\{ |X|^4 \} = \tfrac{8}{(\nu+1)(\nu+3)}$.
\end{definition}

In the Jacobi ensemble, they will be Beta-distributed:

\begin{definition}
A real-valued random variable $X$ is said to be Beta-distributed with parameters $s,t>0$,
which we denote by $X\sim B(s,t)$, if
\begin{equation}\label{E:beta}
\Exp\{f(X)\} = \frac{2^{1-s-t}\Gamma(s+t)}{\Gamma(s)\Gamma(t)}
\int_{-1}^1 f(x) (1-x)^{s-1}(1+x)^{t-1} \, dx.
\end{equation}
Note that $\Exp\{X\}=\frac{t-s}{t+s}$ and $\Exp\{X^2\}=\frac{(t-s)^2+(t+s)}{(t+s)(t+s+1)}$.
\end{definition}

The combined content of Theorem~1.2 and Proposition~B.2 from \cite{KN} is:

\begin{theorem}\label{T:KN_C}
Given $\beta>0$, let $\alpha_k\sim\Theta_{\beta(k+1)+1}$ be independent random variables and let $e^{i\eta}$
be independent and uniformly distributed on $S^1$.  If one forms $B_{n-1}(z)$ according to \eqref{BlaschkeRec},
then
$$
\{ z : B_{n-1}(z) = e^{-i\eta} \}
$$
is distributed according to the C$\beta\!$E${}_n$ ensemble.
\end{theorem}

The result for the Jacobi model is similar; however, we need one extra ingredient.  A probability
measure on $S^1$ is invariant under complex conjugation,
$$
\int f(z) d\mu(z) = \int f(\bar z) d\mu(z),
$$
if and only if it has real Verblunsky coefficients.  From Theorem~1.5 and Proposition~B.2 of \cite{KN}, we have

\begin{theorem}\label{T:KN_J}
Given $\beta>0$, let $\alpha_k$ be independent and distributed as follows
\begin{equation}\label{DistAR}
\alpha_k \sim \begin{cases}
B(\tfrac{k}{4}\beta + a, \tfrac{k}{4}\beta + b)    & \text{$k$ even,} \\
B(\tfrac{k-1}{4}\beta + a + b, \tfrac{k+1}{4}\beta)        & \text{$k$ odd.}
\end{cases}
\end{equation}
Then the $n$ points
\begin{equation}\label{cos of B}
\{ z + \bar z :  B_{2n-1}(z) = -1 \}
\end{equation}
are distributed according to the the Jacobi ensemble \eqref{E:JbetaE}.
\end{theorem}

\subsection{The basic processes}  In the previous subsection, we showed how one could
encode the $\beta$ ensembles as random Blaschke products built from independent random
variables through a simple recurrence, \eqref{BlaschkeRec}.  In this subsection, we will
take matters a few steps further; a similar path was followed in \cite{KStoi}.

From Theorems~\ref{T:KN_C} and~\ref{T:KN_J} we see that
we need only study the Blaschke products on the boundary of the unit circle.  To this end,
let us introduce random continuous functions $\psi_k:(-\pi,\pi)\to\Reals$ via
$B_k(e^{i\theta}) = e^{i\psi_k(\theta)}$.  We break the $2\pi\Ints$ ambiguity by choosing
a branch for the logarithm in the recurrence relation \eqref{BlaschkeRec}.  Specifically,
we define
\begin{equation}\label{psi_rec}
\psi_0(\theta)=\theta, \qquad \psi_{k+1}(\theta) = \psi_{k}(\theta) + \theta + \Upsilon(\psi_k(\theta),\alpha_k) ,
\end{equation}
where 
\begin{equation}\label{Log_series}
\Upsilon(\psi,\alpha) := - 2 \Im \log[ 1-\alpha e^{i\psi} ] =
    \Im \sum_{\ell=1}^\infty \tfrac{2}{\ell} e^{i\ell\psi} \alpha^\ell.
\end{equation}
As the argument of a Blaschke product, $\psi_k$ is an increasing function of $\theta$.

Readers wishing to compare this with \cite{KStoi} should note that here $\psi_k(\theta)$ denotes the
true Pr\"ufer phase rather than the relative Pr\"ufer phase.

To motivate what follows, let us give a quick description of how we will prove Gaussian fluctuations in the circular
case.  The Jacobi case is analogous.

If we choose $\eta,\alpha_k$ to be distributed as in Theorem~\ref{T:KN_C} and use these to produce a sequence of
increasing functions $\psi_k(\theta)$ as in \eqref{psi_rec}, then the set of points
$$
\{ e^{i\theta} : \psi_{n-1}(\theta) \in 2\pi\Ints+\eta \}
$$
will be distributed according to the C$\beta$E$_n$.  In particular, the number of points lying in the arc
$[a,b]\subset(-\pi,\pi)$ is approximately $\frac1{2\pi}[\psi_{n-1}(b)-\psi_{n-1}(a)]$; indeed the error
is plus or minus one.  In this way, it suffices to show that asymptotically, $\psi_{n-1}(b)$ and $\psi_{n-1}(a)$ 
follow a joint normal law.  The error of $\pm1$ will drop out in the limit when we divide by the square root of
the variance.

In the circular case, $\psi_k(\theta)-(k+1)\theta$ is a sum of independent random variables (because $\alpha_k$
have rotationally invariant laws) and so it is easy
to demonstrate that it has Gaussian behaviour.  To study the joint distribution of $\psi_{n-1}$ at several values
of $\theta$ takes a little more care, but can be dealt with using the Central Limit Theorem for Martingales.
Implementing this argument requires a few estimates on $\Upsilon$, which we record here.

\begin{lemma}\label{L:CBasic}
Suppose $\phi,\psi\in\Reals$ and $\alpha\sim\Theta_\nu$.
If $\,\tilde\Upsilon(\psi,\alpha) = 2\Im\{ e^{i\psi} \alpha\}$, then
\begin{gather}
\label{MVT}
     \Exp\{ \Upsilon(\psi, \alpha) \} = \Exp\{ \tilde\Upsilon(\psi,\alpha) \} = 0,\\[2mm]
\label{LogEst2}
    \Exp\bigl\{ \tilde\Upsilon(\psi, \alpha) \tilde\Upsilon(\phi,\alpha) \bigr\} = \tfrac{4}{\nu+1} \cos(\psi-\phi), \\[2mm]
\label{LogEst4}
    \Exp\{ \tilde\Upsilon(\psi,\alpha)^4 \} = \tfrac{48}{(\nu+1)(\nu+3)}, \\[2mm]
\label{C_diff}
    \Exp\bigl\{ \bigl| \Upsilon(\psi, \alpha) - \tilde\Upsilon(\psi,\alpha) \bigr|^2 \bigr\}
        \leq \tfrac{16}{(\nu+1)(\nu+3)}.
\end{gather}
Lastly, combining \eqref{LogEst2} and \eqref{C_diff} shows
\begin{gather}\label{LogEst2b}
\Exp\bigl\{ \bigl| \Upsilon(\psi, \alpha) \bigr|^2 \bigr\} \leq \tfrac{8}{(\nu+1)}.
\end{gather}
\end{lemma}

\begin{proof}
The vanishing of $\Exp\{ \Upsilon(\psi, \alpha) \}$ merely uses the fact that $\alpha$ follows a rotationally invariant law.
Specifically, for any $0\leq r < 1$,
\begin{gather*}
 \int_0^{2\pi} \log\bigl[ 1 - re^{i\theta+i\psi} \bigr] \, \tfrac{d\theta}{2\pi} = 0
\end{gather*}
by the Mean Value Principle for harmonic functions.  The second identity in \eqref{MVT} follows in the same manner.

To prove \eqref{LogEst2} we simply compute:
\begin{align*}
\Exp\bigl\{ \tilde\Upsilon(\psi, \alpha) \tilde\Upsilon(\phi,\alpha) \bigr\}
    &=\tfrac{4(\nu-1)}{2\pi} \!\int\!\!\!\int_\Disk r^2\sin(\theta+\psi)\sin(\theta+\phi) (1-r^2)^{(\nu-3)/2} r\,dr\,d\theta \\
    &= \tfrac{4}{\nu+1}\cos(\psi-\phi).
\end{align*}
Similarly,
\begin{align}
\Exp\bigl\{ \tilde\Upsilon(\psi, \alpha)^4 \bigr\}
    &= \tfrac{16(\nu-1)}{2\pi} \!\int\!\!\!\int_\Disk r^4\sin^4(\theta+\psi) (1-r^2)^{(\nu-3)/2} r\,dr\,d\theta \\
    &= \tfrac{48}{(\nu+1)(\nu+3)}.
\end{align}

Applying Plancharel's theorem to the power series formula for $\Upsilon$ gives
\begin{align}
\Exp\bigl\{ \bigl| \Upsilon(\psi, \alpha) - \tilde\Upsilon(\psi,\alpha) \bigr|^2 \bigr\}
    &= \sum_{\ell=2}^\infty \tfrac{2}{\ell^2} \Exp\{ |\alpha|^{2\ell} \}
    \leq 2(\tfrac{\pi^2}{6}-1) \Exp\{ |\alpha|^{4} \}.
\end{align}
Noting that $\tfrac{\pi^2}{6}\leq 2$ and  $\Exp\{ |\alpha|^{4} \} = \tfrac{8}{(\nu+1)(\nu+3)}$ proves \eqref{C_diff}.
\end{proof}

The analogue of Lemma~\ref{L:CBasic} for the Jacobi case is

\begin{lemma}\label{L:JBasic}
Let $\tilde \Upsilon(\psi,\alpha)=2[\alpha-\Exp(\alpha)]\sin(\psi)$ with $\alpha\sim B(s,t)$.  Then
\begin{align}
\label{J_MV}
    \Exp\{ \tilde \Upsilon(\psi, \alpha) \} &= 0\\
\label{JEst2}
    \Exp\bigl\{ \tilde \Upsilon(\psi, \alpha) \tilde \Upsilon(\phi,\alpha) \bigr\}
        &= \tfrac{8st}{(s+t)^2(s+t+1)} [\cos(\psi-\phi)-\cos(\psi+\phi)]
\end{align}
for any $\phi,\psi\in\Reals$.  Moreover, if $\alpha_k$ is distributed as in \eqref{DistAR}, then
\begin{align}
\label{JEst4}
    \Exp\{ \tilde\Upsilon(\psi,\alpha_k)^4 \} 
    &\lesssim (k+1)^{-2} \\
\label{JEst5}
    \Exp\{ |\Upsilon(\psi,\alpha_k)-2\alpha_k\sin(\psi)&-\alpha_k^2\sin(2\psi) | \} 
    \lesssim (k+1)^{-2} \\
\label{JEst1}    
    \Exp\{ |\Upsilon(\psi, \alpha_k)| \} &\lesssim (k+1)^{-1/2}
\end{align}
where the implicit constant does not depend on $k$, merely $\beta$, $a$, and $b$.
\end{lemma}

\begin{proof}
Equation \eqref{J_MV} is immediate from the definition, while \eqref{JEst2} follows from
\eqref{AVar} and $2\sin(\psi)\sin(\phi) = \cos(\psi-\phi) - \cos(\psi-\phi)$.

To obtain \eqref{JEst4}, we can estimate rather crudely:
$$
\Exp\{ \tilde\Upsilon(\psi,\alpha_k)^4 \} \leq 2^8 \Exp\{ \alpha_k^4 \} \lesssim (k^2+1)^{-1}
$$
using \eqref{A4}.  Next we prove \eqref{JEst5}.  It follows from combining
\begin{align*}
\text{LHS\eqref{JEst5}} = \Exp\biggl\{ \biggl| \sum_{l=3}^\infty \tfrac{2}{l} \sin(l\psi) \alpha^l \biggr| \biggr\}
&\leq \bigl| \Exp\{ \alpha^3 \} \bigr| + \Exp\biggl\{ \sum_{m=2}^\infty [\tfrac{2}{2m} + \tfrac{2}{2m+1} ] \alpha^{2m} \biggr\} \\
&\leq \bigl| \Exp\{ \alpha^3 \} \bigr| + 2\Exp\biggl\{ -\alpha^{2}\log[1-\alpha^{2}] \biggr\}.
\end{align*}
with \eqref{A3} and Lemma~\ref{L:log}.

The last estimate follows from \eqref{JEst5} and \eqref{A2} by applying the triangle and Cauchy--Schwarz inequalities.
\end{proof}

On several occasions we will make use of

\begin{lemma}\label{L:Sum}
Given real valued sequences $\epsilon_k,X_k,Y_k$ with $X_{k+1} = X_k + \delta + Y_k$ for some $\delta\in(0,2\pi)$,
\begin{equation}\label{E:Sum}
\biggl| \sum \epsilon_k e^{i X_k} \biggr| \leq \frac{2\|\epsilon_k\|_{\ell^\infty} +
        \|\epsilon_k - \epsilon_{k-1}\|_{\ell^1} + \|\epsilon_k Y_k \|_{\ell^1}}{|1 - e^{i\delta}|}.
\end{equation}
\end{lemma}

\begin{proof}
Using summation by parts and other elementary manipulations,
\begin{align*}
\sum_{k=0}^{n-1} \epsilon_k e^{i X_k}
    &= \sum_{k=0}^{n-1} \epsilon_k \frac{e^{i X_k}-e^{i (X_k+\delta)}}{1-e^{i\delta}} \\
&= \sum_{k=0}^{n-1} \epsilon_k \frac{e^{i X_k}-e^{i X_{k+1}}}{1-e^{i\delta}} +
    \sum_{k=0}^{n-1} \epsilon_k \frac{e^{i X_{k+1}}-e^{i (X_k+\delta)}}{1-e^{i\delta}}  \\
&= \frac{\epsilon_ne^{iX_n} - \epsilon_0e^{iX_0}}{1-e^{i\delta}} + \sum_{k=1}^{n-1} \frac{\epsilon_k-\epsilon_{k-1}}{1-e^{i\delta}}
    + \sum_{k=0}^{n-1} \epsilon_k \frac{e^{i (X_k+\delta)}[e^{iY_k}-1]}{1-e^{i\delta}} .
\end{align*}
This implies \eqref{E:Sum} in a most obvious way.
\end{proof}

\subsection{Martingale Central Limit Theorem}
In both the Circular and Jacobi cases, we will approximate $\psi_{m-1}(\theta)$ by
\begin{equation}\label{SDefn}
S(m,\theta) = \sum_{k=0}^{m-1} \tilde \Upsilon(\psi_k(\theta),\alpha_k).
\end{equation}
The recursive definition of $\psi_k$ and the fact that $\Exp\{\tilde\Upsilon\}=0$ shows that this is a
martingale for the sigma algebras
$$
\mathcal{M}_m := \sigma(\alpha_0,\ldots,\alpha_{m-1}),
$$
that is, the smallest sigma algebra generated by these random variables.  In the circular case, there
is no need to pass from $\Upsilon$ to $\tilde\Upsilon$ in order to obtain a martingale; nevertheless,
doing so simplifies the proof just a little.

\begin{prop}\label{MCLT}
Fix $\theta_1,\ldots,\theta_J$ distinct. Suppose that as $n\to\infty$
\begin{gather}
\label{MStability}
\frac{1}{\log(n)} \sum_{k=0}^{n-1} \Exp\bigl\{ \tilde\Upsilon(\psi_k(\theta_j),\alpha_k)
    \tilde\Upsilon(\psi_k(\theta_l),\alpha_k) \big| \mathcal{M}_k \bigr\}
    \to \sigma^2 \delta_{jl} \quad \text{in $L^1$ sense} \\
\label{MMoment}
\frac{1}{\log^2(n)} \sum_{k=0}^{n-1} \Exp\{ |\tilde\Upsilon(\psi_k(\theta_j),\alpha_k)|^4 \}
    \to 0  \\
\label{MApprox}
\frac{1}{\sqrt{\log(n)}} \Exp\biggl\{ \biggl| \sum_k \Upsilon(\psi_k(\theta_j),\alpha_k)
    - \tilde\Upsilon(\psi_k(\theta_j),\alpha_k) \biggl| \biggl\}
    \to 0
\end{gather}
for all $1\leq j,l \leq J$.  Then in the random variables
\begin{equation}\label{Mconclusion}
\frac{\psi_{n-1}(\theta_j) - n\theta_j}{\sqrt{\log(n)}}, \qquad j\in\{1,2,\ldots,J\},
\end{equation}
converge to independent Gaussians of mean zero and variance $\sigma^2$.
\end{prop}

\begin{proof}
The first two conditions are far more than is really needed to apply the Martingale Central Limit Theorem
to the sequence of random vectors
$$
\tfrac{1}{\sqrt{\log(n)}} \bigl( S(n,\theta_1),\ldots,S(n,\theta_J)\bigr);
$$
see \cite{Pollard,Varadhan}, for example.  These references present the proof for scalar martingales --- a small
extension of the usual Central Limit Theorem --- which can then be applied to any linear combination of the components
of the vector.  In this way, we see that
$$
\tfrac{1}{\sqrt{\log(n)}} \bigl( S(n,\theta_1),\ldots,S(n,\theta_J)\bigr) \to N(0,\sigma^2\delta_{jl}) 
$$
in distribution. This extends to \eqref{Mconclusion} because \eqref{MApprox} says that the difference
between these to vectors converges to $0$ in $L^1$ sense.
\end{proof}

\section{Circular case}

In this section, we proof Theorem~\ref{T:mainC}.  The discussion in the previous section explains why
this amounts to checking the hypotheses in Proposition~\ref{MCLT}.   We begin with a recap.

Let $\alpha_k\sim\Theta_{\beta(k+1)+1}$ be chosen independently (as in Theorem~\ref{T:KN_C}) and then
form the process $\psi_{n-1}(\theta)$ by solving the recurrence \eqref{psi_rec}.  If $\eta$ is chosen
independently from $[0,2\pi)$ according to the uniform distribuion, then
$$
\{ e^{i\phi} : \psi_{n-1}(\phi)+\eta \in 2\pi\Ints \} 
$$
is distributed according to C$\beta$E$_n$.  Notice that the number of points in the arc $[a,b]\subset(-\pi,\pi)$
differs from $\tfrac{1}{2\pi}[\psi_{n-1}(b)-\psi_{n-1}(a)]$ by at most $\pm 1$.  Thus Theorem~\ref{T:mainC}
will follow once we show that
$$
\tfrac{1}{\sqrt{\log(n)}}[\psi_{n-1}(\theta)-n\theta] \to \sqrt{\tfrac{4}{\beta}}\,\Psi(\theta)
$$
in distribution.  This in turn follows from Proposition~\ref{MCLT} once we check its hypotheses.
The two lemmas that follow verify these conditions in the order they appear there. 

\begin{lemma}\label{L:StabC}  Given $\theta_1,\theta_2\in(-\pi,\pi)$,
\begin{gather}
\label{MStabC}
\frac{1}{\log(n)} \sum_{k=0}^{n-1} \Exp\bigl\{ \tilde\Upsilon(\psi_k(\theta_1),\alpha_k)
    \tilde\Upsilon(\psi_k(\theta_2),\alpha_k) \big| \mathcal{M}_k \bigr\}
    \to \begin{cases} \tfrac{4}{\beta} &: \theta_1=\theta_2 \\ 0 &: \theta_1\neq\theta_2 \end{cases}
\end{gather}
in $L^1$ sense.
\end{lemma}

\begin{proof}
By \eqref{LogEst2},
\begin{align}\label{MStabC2}
\text{LHS\eqref{MStabC}} = \frac{1}{\log(n)} \sum_{k=0}^{n-1} 
\tfrac{4}{\beta(k+1)+2} \cos\bigl(\psi_k(\theta_1)-\psi_k(\theta_2)\bigr),
\end{align}
which immediately settles the case $\theta_1=\theta_2$. By symmetry, this leaves us only to treat the case
$\theta_1 > \theta_2$.  We do this using Lemma~\ref{L:Sum} with
$$
X_k=\psi_k(\theta_1)-\psi_k(\theta_2), \quad
Y_k=\Upsilon(\psi_k(\theta_1),\alpha_k) - \Upsilon(\psi_k(\theta_2),\alpha_k), \quad
\epsilon_k = \tfrac{4}{\beta(k+1)+2},
$$
and $\delta=\theta_1-\theta_2$.  Indeed, we may deduce
$$
\Exp\Bigl\{ \Bigl| \text{ RHS\eqref{MStabC2} } \Bigr| \Bigr\} = O(1/\log(n))
$$
since by \eqref{LogEst2}, $\Exp\{ |Y_k| \} \leq 8[\beta(k+1)+2)]^{-1/2}$.
\end{proof}

\begin{lemma}  Given $\theta\in(-\pi,\pi)$,
\begin{gather}\label{MMomC}
\frac{1}{\log^2(n)} \sum_{k=0}^{n-1} \Exp\{ |\tilde\Upsilon(\psi_k(\theta),\alpha_k)|^4 \}
    \to 0\\
\label{MApprC}
\frac{1}{\log(n)} \Exp\biggl\{ \biggl| \sum_{k=0}^{n-1} \Upsilon\bigl(\psi_k(\theta),\alpha_k\bigr)
    - \tilde\Upsilon\bigl(\psi_k(\theta),\alpha_k\bigr) \biggl|^2 \biggl\}
    \to 0 
\end{gather}
as $n\to\infty$.
\end{lemma}

\begin{proof}
From \eqref{LogEst4}, we see that not only does LHS\eqref{MMomC} converge to zero, it is $O(1/\log^2(n))$.
We turn now to \eqref{MApprC}.  By \eqref{MVT} and \eqref{C_diff} we have
\begin{align*}
\text{LHS\eqref{MApprC}} &= \frac{1}{\log(n)} \sum_{k=0}^{n-1}  \Exp\bigl\{ | \Upsilon(\psi_k(\theta_j),\alpha_k)
    - \tilde\Upsilon(\psi_k(\theta_j),\alpha_k) |^2 \bigl\} \\
&\leq \frac{1}{\log(n)} \sum_{k=0}^{n-1} \frac{16}{[\beta(k+1)+2][\beta(k+1)+4]},
\end{align*}
which is $O(1/\log(n))$.  Note that \eqref{MApprC} implies \eqref{MApprox} via the Cauchy--Schwarz inequality.
\end{proof}

\section{Jacobi case}

In the Jacobi case we choose $\alpha_k$ to be distributed as in Theorem~\ref{T:KN_J}.
For then
$$
\bigl\{ 2\cos(\theta) : \psi_{n-1}(\theta) - \pi \in 2\pi\Ints,\ \theta\in(0,\pi) \bigr\}
$$
is distributed according to J$\beta$E$_n$ with parameters $a$ and $b$.  As the Verblunsky
coefficients are real-valued, $\psi_k(0)\equiv 0$.  Thus 
$$
| N_n(\theta) - \tfrac1{2\pi}\psi_{n-1}(\theta) | \leq \tfrac12 
$$
with probability one. Recall from the introduction that $N_n(\theta)$ denotes the number of particles
in the interval $[2\cos(\theta),2]$.

In light of this discussion, we see that in order to prove Theorem~\ref{T:mainJ}, we need only show
that for any distinct $\theta_1,\ldots,\theta_J\in(0,\pi)$,
$$
\sqrt{\tfrac{\beta}{4\log(n)}}\,\bigl[\psi_{n-1}(\theta_j) - n\theta_j\bigr] 
$$
converge to independent Gaussian random variables with mean $0$ and variance $1$ as $n\to\infty$.
This in turn can be effected via Proposition~\ref{MCLT}
provided we verify its hypotheses.  That is precisely what we will do.

\begin{lemma}  Given $\theta_1,\theta_2\in(0,\pi)$,
\begin{gather}
\label{MStabJ}
\frac{1}{\log(n)} \sum_{k=0}^{n-1} \Exp\bigl\{ \tilde\Upsilon(\psi_k(\theta_1),\alpha_k)
    \tilde\Upsilon(\psi_k(\theta_2),\alpha_k) \big| \mathcal{M}_k \bigr\}
    \to \begin{cases} \tfrac{4}{\beta} &: \theta_1=\theta_2 \\ 0 &: \theta_1\neq\theta_2 \end{cases}
\end{gather}
in $L^1$ sense.
\end{lemma}

\begin{proof}
By \eqref{JEst2},
\begin{align}\label{MStabJ2}
\text{LHS\eqref{MStabJ}} = \frac{1}{\log(n)} \sum_{k=0}^{n-1} 
    \epsilon_k \bigl[\cos\bigl(\psi_k(\theta_1)-\psi_k(\theta_2)\bigr)-\cos\bigl(\psi_k(\theta_1)+\psi_k(\theta_2)\bigr)\bigr],
\end{align}
where
$$
\epsilon_k = \begin{cases} \frac{4(k\beta+4a)(k\beta+4b)}{(k\beta+2a+2b)^2(k\beta+2a+2b+2)} &: k \text{ even} \\[2ex]
    \frac{4[(k-1)\beta+4a+4b](k+1)\beta}{(k\beta+2a+2b)^2(k\beta+2a+2b+2)} &: k \text{ odd.}\end{cases}
$$

Applying Lemma~\ref{L:Sum} with 
$$
X_k=\psi_k(\theta_1)+\psi_k(\theta_2), \quad
Y_k=\Upsilon(\psi_k(\theta_1),\alpha_k) + \Upsilon(\psi_k(\theta_2),\alpha_k), \quad
$$
and $\delta=\theta_1+\theta_2$, then using \eqref{JEst1}, shows that
\begin{align}\label{MStabJ23}
\text{LHS\eqref{MStabJ}} = \frac{1}{\log(n)} \sum_{k=0}^{n-1} 
    \epsilon_k \cos\bigl(\psi_k(\theta_1)-\psi_k(\theta_2)\bigr) + O\bigl(1/\log(n)\bigr).
\end{align}
This renders \eqref{MStabJ} in the same manner as in Lemma~\ref{L:StabC}.
\end{proof}

\begin{lemma}  Given $\theta\in(0,\pi)$,
\begin{gather}\label{MMomJ}
\frac{1}{\log^2(n)} \sum_{k=0}^{n-1} \Exp\{ |\tilde\Upsilon(\psi_k(\theta),\alpha_k)|^4 \}
    \to 0\\
\label{MApprJ}
\frac{1}{\sqrt{\log(n)}} \Exp\biggl\{ \biggl| \sum_{k=0}^{n-1} \Upsilon\bigl(\psi_k(\theta),\alpha_k\bigr)
    - \tilde\Upsilon\bigl(\psi_k(\theta),\alpha_k\bigr) \biggl| \biggl\}
    \to 0 
\end{gather}
as $n\to\infty$.
\end{lemma}

\begin{proof}
By \eqref{JEst4}, not only does LHS\eqref{MMomJ} converge to zero, it is $O(1/\log^2(n))$.

We turn now to \eqref{MApprJ}.  By \eqref{JEst5},
\begin{align*}
\text{LHS\eqref{MApprJ}} &= \frac{1}{\sqrt{\log(n)}} \Exp\biggl\{ \biggl| \sum_{k=0}^{n-1} 2\Exp\{\alpha_k\}\sin(\psi_k(\theta))
    - \alpha_k^2\sin(2\psi_k(\theta)) \biggl| \biggl\}
    + O\Bigl(\tfrac1{\sqrt{\log(n)}}\Bigr).
\end{align*}
Lemma~\ref{L:Sum} with $\delta=\theta$, $X_k=\psi_k(\theta)$, and $\epsilon_k=\Exp\{\alpha_k\}$ shows that
\begin{align*}
\Exp\biggl\{ \biggl| \sum_{k=0}^{n-1} 2\Exp\{\alpha_k\}\sin(\psi_k(\theta)) \biggr| \biggr\} =O(1)
\end{align*}
and in a similar way,
\begin{align*}
\Exp\biggl\{ \biggl| \sum_{k=0}^{n-1} \Exp\{\alpha_k^2\}\sin(2\psi_k(\theta)) \biggr| \biggr\} = O(1).
\end{align*}
In view of these estimates, \eqref{MApprJ} follows from
\begin{align*}
\Exp\biggl\{ \biggl| \sum_{k=0}^{n-1} \bigl[ \alpha_k^2 - \Exp\{\alpha_k^2\} \bigr] \sin(2\psi_k(\theta)) \biggl|^2 \biggl\}
&=\sum_{k=0}^{n-1} \Exp\Bigl\{ \bigl[ \alpha_k^2 - \Exp\{\alpha_k^2\} \bigr]^2 \sin^2(2\psi_k(\theta)) \Bigr\} \\
&\leq\sum_{k=0}^{n-1} \Exp\bigl\{ \alpha_k^4 \bigr\} = O(1).
\end{align*}
The last deduction is based on \eqref{A4}.
\end{proof}

\appendix
\section{Some integrals of Beta random variables}

As seen in \eqref{E:beta} above,
$$
I(s,t):= \int_{-1}^1 (1-x)^{s-1}(1+x)^{t-1} \, dx = \frac{2^{s+t-1}\Gamma(s)\Gamma(t)}{\Gamma(s+t)};
$$
indeed this is the famous Beta integral of Euler.  In this appendix, we record a few simple computations
that were needed in the text; nothing is novel.

By the binomial theorem,
$$
x^k = 2^{-k} \sum_{m=0}^k (-1)^m \binom{k}{m}  (1-x)^{m} (1+x)^{k-m} 
$$
and so if $X\sim B(s,t)$, then
$$
\Exp\{ X^k \} = \sum_{m=0}^k (-1)^m \binom{k}{m} \frac{I(s+m,t+k-m)}{I(s,t)}.
$$
With a little algebra one then obtains
\begin{align}
\label{A1} \Exp\{X\}   &=\frac{t-s}{t+s} \\
\label{A2} \Exp\{X^2\} &=\frac{(t-s)^2+(t+s)}{(t+s)(t+s+1)} \\
\label{A3} \Exp\{X^3\} &=\frac{(t-s)[(t-s)^2+3(t+s)+2]}{(s+t)(1+s+t)(s+t+2)}\\
\label{A4} \Exp\{X^4\} &=\frac{(t-s)^2[(t-s)^2+6(s+t)+8] +3(t+s)^2+6(t+s)}{(s+t)(1+s+t)(s+t+2)(s+t+3)}.
\end{align}
In particular,
\begin{align}
\label{AVar} \Exp\bigl\{ \bigl( X -\Exp\{X\} \bigr)^2 \bigr\} &=\frac{4st}{(t+s)^2(t+s+1)}.
\end{align}

Differentiating $I(s,t)$ with respect to $s$ or $t$ gives access to expectations of logarithms.  For example,
\begin{align*}
\Exp\{ - X^2\log[1-X^2] \} &= \frac{(\partial_s+\partial_t) [I(s+1,t+1) + I(s,t)]}{I(s,t)} \\
&= \frac{(s-t)^2+s+t}{(s+t)(1+s+t)}[2\Psi(s+t)-\Psi(t)-\Psi(s)-2\log(2)] \\
& \qquad\qquad         + 4\frac{(s+t)(s-t)^2+t^2+s^2}{(s+t)^2(1+s+t)^2}.
\end{align*}
Where $\Psi(x):=\partial_x \log\circ\,\Gamma(x)$ is the digamma function.

\begin{lemma}\label{L:log}
Given $s,t\in(\epsilon,\infty)$ obeying $|s-t|\leq \delta$ and $X\sim B(s,t)$,
\begin{align}
\Exp\{ - X^2\log[1-X^2] \} &\lesssim  (s+t)^{-2}
\end{align}
where the implicit constant depends only on $\epsilon$ and $\delta$.
\end{lemma}

\begin{proof}
In light of the calculation above, we need only prove
$$
|2\Psi(s+t)-\Psi(t)-\Psi(s)-2\log(2)| \lesssim (s+t)^{-1}.
$$
This in turn can be deduced from $\Psi(x)=\log(x)+O(x^{-1})$ as $x\to\infty$ --- a very weak form
of Stirlings formula.  Indeed,
$$
2\log(s+t)-\log(t)-\log(s)-2\log(2) = \log\left[\tfrac{(s+t)^2}{4st}\right] =  \log\left[ 1 + \tfrac{(s-t)^2}{4st}\right]
$$ 
which is $\lesssim (s+t)^{-2}$. 
\end{proof}



\begin{thebibliography}{10}
\newcommand{\MSN}[1]{\href{http://www.ams.org/mathscinet-getitem?mr=#1}{\sc MR#1}}

\bibitem{CMV} M.~J.~Cantero, L.~Moral, and L.~Vel\'azquez,
Five-diagonal matrices and zeros of orthogonal polynomials on the unit circle. \textit{Linear
Algebra Appl.} {\bf 362} (2003), 29--56.
\MSN{1955452}

\bibitem{CostLeb} O.~Costin and J. L.~Lebowitz,
Gaussian fluctuation in random matrices.
\textit{Phys. Rev. Lett.} \textbf{75} (1995), 69--72.




\bibitem{DiaEvans} P.~Diaconis and S.~Evans,
Linear Functionals of Eigenvalues of Random Matrices.
\textit{Trans. AMS} \textbf{353} (2001) 2615--2633.
\MSN{1828463}

\bibitem{DiaShah} P.~Diaconis and M.~Shahshahani,
On the Eigenvalues of Random Matrices,
\textit{Journal of Applied Probability} \textbf{31} (1994) 49--61
\MSN{1274717}

\bibitem{DumE} I.~Dumitriu and A.~Edelman,
Matrix models for beta ensembles.
\textit{J. Math. Phys.} {\bf 43} (2002), 5830--5847.
\MSN{1936554}

\bibitem{DumE2} I.~Dumitriu and A.~Edelman,
Global Spectrum Fluctuations for the $\beta$-Hermite and $\beta$-Laguerre ensembles via matrix models.
\textit{J. Math. Phys.} {\bf 47} (2006).
\MSN{2239975}

\bibitem{Dyson} F.~Dyson,
Statistical theory of the energy levels of complex systems. I, II, and III. \textit{J. Math. Phys.}
{\bf 3} (1962), 140--156, 157--165, and 166--175. \MSN{0143556}, \MSN{0143557}, \MSN{0143558}.

\bibitem{ForresterBook} P.~J.~Forrester,
\textit{Log-gases and Random matrices.}  Available from the author's web page.




\bibitem{Good} I.~J.~Good,
Short proof of a  conjecture by Dyson. \textit{J. Math. Phys.} {\bf 11} (1970), 1884.
\MSN{0258644}


\bibitem{Joh1} K.~Johannson,
On Szeg\H{o}'s asymptotic formula for Toeplitz determinants and generalizations.
\textit{Bull. Sci. Math} \textbf{112} (1988), 257--304.
\MSN{0975365}

\bibitem{Joh2} K.~Johannson,
On random matrices from the classical compact groups.
\textit{Ann. of Math.} \textbf{145} (1997), 519--545.
\MSN{1454702}

\bibitem{Joh3} K.~Johannson,
On fluctuations of eigenvalues of random Hermitian matrices.
\textit{Duke Math. J.} \textbf{91} (1998), 151--204.
\MSN{1487983}


\bibitem{Jonsson} D.~Jonsson,
Some limit theorems for the eigenvalues of a sample covariance matrix.
\textit{J. Multivariate Anal.} \textbf{12} (1982), 1--38. 
\MSN{0650926}

\bibitem{KN} R.~Killip and I.~Nenciu,
Matrix models for circular ensembles. \textit{Int. Math. Res. Not.} \textbf{2004}, 2665--2701.
\MSN{2127367}


\bibitem{KStoi} R.~Killip and M.~Stoiciu,
Eigenvalue Statistics for CMV Matrices: from Poisson to Clock via C$\beta$E.  Preprint \texttt{math-ph/0608002}.




\bibitem{Mehta} M.~L.~Mehta,
\textit{Random matrices.} Third Edition. Pure and Applied Mathematics (Amsterdam), 142.
Elsevier/Academic Press, Amsterdam, 2004.



\bibitem{Pollard} D.~Pollard, \emph{Convergence of stochastic processes.}
Springer Series in Statistics. Springer-Verlag, New York, 1984.
\MSN{0762984}

\bibitem{RRV} J. Ramirez, B. Rider, and B. Vir\'ag,
Beta ensembles, stochastic Airy spectrum, and a diffusion.
Preprint \texttt{math.PR/0607331}.

\bibitem{Selberg} A.~Selberg,
Bemerkninger om et multipelt integral.
\textit{Norsk Mat. Tidsskr.} {\bf 26} (1944), 71--78.
\MSN{0018287}


\bibitem{SimonOPUC1} B.~Simon,
\textit{Orthogonal Polynomials on the Unit Circle, vol. 1}. American Mathematical Society
Colloquium Publications, American Mathematical Society, Providence, Rhode Island, 2004.

\bibitem{SimonOPUC2} B.~Simon,
\textit{Orthogonal Polynomials on the Unit Circle, vol. 2}. American Mathematical Society
Colloquium Publications, American Mathematical Society, Providence, Rhode Island, 2004.




\bibitem{Sosh1} A.~Soshnikov,
Determinantal Random Point Fields.
\textit{Russian Mathematical Surveys} \textbf{55} (2000), 923--975.
\MSN{1799012}

\bibitem{Sosh2} A.~Soshnikov,
Gaussian Fluctuation for the Number of Particles in Airy, Bessel, Sine and Other Determinantal Random Point Fields.
\textit{Journal of Stat. Phys.} \textbf{100} (2000), 491--522.
\MSN{1788476}

\bibitem{Sosh3} A.~Soshnikov,
Gaussian limit for determinantal random point fields.
\textit{Ann. Probab.} \textbf{30} (2002), 171--187.
\MSN{1894104} 


\bibitem{SzegoToe} G.~Szeg\H{o}, 
On certain Hermitian forms associated with the Fourier series of a positive function.
\textit{Comm. S\'em. Math. Univ. Lund}  (1952). Tome Supplementaire, 228--238.
\MSN{0051961}

\bibitem{Szego} G.~Szeg\H{o},
\textit{Orthogonal Polynomials.} American Mathematical Society Colloquium Publications, Vol. XXIII.
American Mathematical Society, Providence, Rhode Island, 1975.

\bibitem{Trotter} H.~F.~Trotter,
Eigenvalue distributions of large Hermitian matrices; Wigner's semicircle law and a theorem of Kac, Murdock, and Szeg\H{o}.
\textit{Adv. in Math.} \textbf{54} (1984), 67--82.
\MSN{0761763}

\bibitem{Varadhan} S. R. S. Varadhan, \textit{Probability theory.} Courant Lecture Notes in Mathematics, \textbf{7}.
New York University, Courant Institute of Mathematical Sciences, New York; American Mathematical Society, Providence, RI.
\MSN{1852999}

\bibitem{VV} B. Valk\'o and B. Vir\'ag, \textit{Scaling limits of random matrices: the stochastic sine equation.}
In preparation.



\bibitem{Wieand}  K.~Wieand,
Eigenvalue distributions of random unitary matrices.
\textit{Probab. Theory Related Fields} \textbf{123} (2002), 202--224. 
\MSN{1900322}

\bibitem{Wilson} K.~Wilson,
Proof of a conjecture by Dyson. \textit{J. Math. Phys.} {\bf 3} (1962), 1040--1043.
\MSN{0144627}

\end{thebibliography}
\end{document}